\documentclass[conference,letterpaper]{IEEEtran}

\usepackage[utf8]{inputenc} 
\usepackage[T1]{fontenc}
\usepackage{url}
\usepackage{ifthen}
\usepackage{cite}
\usepackage[cmex10]{amsmath} 
\usepackage{amssymb}

\newcommand{\cA}{{\mathcal{A}}}

\newcommand{\cL}{{\mathcal{L}}}

\newcommand{\cP}{{\mathcal{P}}}

\newcommand{\cU}{{\mathcal{U}}}

\newcommand{\cX}{{\mathcal{X}}}

\newcommand{\cY}{{\mathcal{Y}}}

\usepackage[mathcal]{euscript}

\DeclareMathOperator{\ev}{\mathbb{E}}
\DeclareMathOperator{\pr}{\mathbb{P}}

\newcounter{actr}
{\begin{list}{(\alph{actr})}{\usecounter{actr}}}{\end{list}}

\newcounter{ictr}
{\begin{list}{(\roman{ictr})}{\usecounter{ictr}}}{\end{list}}

{

}%
{

}

\newcommand{\abs}[1]{\left|#1\right|}

\newcommand{\defeq}{\triangleq}

\newcommand{\E}[1]{\ev\left[{#1}\right]}
\newcommand{\Ed}[2]{\ev_{#1}\left[{#2}\right]} 

\newcommand{\Prob}[1]{\pr\left({#1}\right)}
\newcommand{\prob}[1]{\Prob{#1}}

\newcommand{\reals}{\mathbb{R}}

\newcommand{\thmref}[1]{Theorem~\ref{#1}}
\newcommand{\corolref}[1]{Corollary~\ref{#1}}

\newcommand{\lemref}[1]{Lemma~\ref{#1}}

\usepackage{todonotes}

\newcommand{\loss}{\mathcal{L}}
\newcommand{\encampsParant}[1]{\left(#1\right)}
\newcommand{\littleOh}[1]{o\encampsParant{#1}}
\newcommand{\bigOh}[1]{O\encampsParant{#1}}

\DeclareMathOperator{\bin}{Bin}
\newcommand{\cardx}{\abs{\cX}}
\newcommand{\cardy}{\abs{\cY}}
\newcommand{\cq}{\check{q}}

\newtheorem{thm}{Theorem}
\newtheorem{lemma}{Lemma}

\newtheorem{corol}{Corollary}

\allowdisplaybreaks

\begin{document}
\title{On Semi-supervised Estimation  \\ of Discrete Distributions
under $f$-divergences 
} 


 \author{%
  \IEEEauthorblockN{H.S.Melihcan Erol and Lizhong Zheng}
  \IEEEauthorblockA{Dept. EECS and RLE\\
                    Massachusetts Institute of Technology\\
                    Cambridge, MA 02139\\
                    \{hsmerol, lizhong\}@mit.edu}
 }



\maketitle


\begin{abstract}
  We study the problem of estimating the joint probability mass function (pmf) over two random variables. In particular, the estimation is based on the observation of $m$ samples containing both variables and $n$ samples missing one fixed variable. We adopt the minimax framework with $l^p_p$ loss functions.  Recent work established that univariate minimax estimator combinations achieve minimax risk with the optimal first-order constant for $p \ge 2$ in the regime $m = o(n)$, questions remained for $p \le 2$ and various $f$-divergences. In our study, we affirm that these composite estimators are indeed minimax optimal for $l^p_p$ loss functions, specifically for the range $1 \le p \le 2$, including the critical $l_1$ loss. Additionally, we ascertain their optimality for a suite of $f$-divergences, such as KL, $\chi^2$, Squared Hellinger, and Le Cam divergences.
\end{abstract}

\section{Introduction}
Estimating the probability mass function (pmf) is a pivotal task in statistics, and the minimax methodology has emerged as a fundamental approach in this context, a fact well-established in \cite{Wald49}. Pioneering studies adeptly navigated the realm of minimax risk in pmf estimation under $l^2_2$ loss, culminating in the identification of the quintessential minimax estimator \cite{trybula1958some,Olkin1979,wilczynski1985minimax}. Further research refined these insights by accurately determining the primary constant for the minimax risk concerning KL-divergence, $l_1$, and $f$-divergences \cite{braess2004bernstein, HanJW14,pmlr-v40-Kamath15}.

The recent growth in dataset sizes often exceeds the capacity for thorough labeling, leading to datasets with a significant number of samples missing specific information. The estimators discussed in studies such as \cite{trybula1958some, Olkin1979, wilczynski1985minimax, braess2004bernstein, HanJW14, pmlr-v40-Kamath15} fall short as they fail to effectively leverage these unlabeled samples, resulting in suboptimal performance. This challenge underscores the importance of refined semi-supervised learning techniques that adeptly combine both labeled and unlabeled data. This paper strives to contribute to this area by enhancing the approach introduced in \cite{onsemsup}, extending its reach to include a spectrum of $f$-divergences including KL, $\chi^2$, Squared Hellinger, and Le Cam divergences.

Echoing the framework in \cite{onsemsup}, our analysis is concentrated on two dependent random variables, $X$ and $Y$, with a distribution denoted by $p_{XY}$. We scrutinize two distinct datasets: the first composed of $m$ i.i.d. samples of $(x_i, y_i)$ sourced from $p_{XY}$, and the second consisting of $n$ samples of $x_j$ derived from the marginal distribution $p_X$. This work extends and enriches the foundational principles set forth in \cite{onsemsup} by demonstrating that the composition of univariate minimax estimators retains the first-order risk term's accuracy across $l^p_p$ loss functions for the range $1 \le p \le 2$ and encompasses a broad array of $f$-divergences. This spectrum includes, notably, total-variation, $\chi^2$, Squared Hellinger, and Lecam divergences, thereby broadening the scope of the initial framework which focused on $l^p_p$ loss functions for $p \ge 2$.
\section{Notations \& Preliminaries}
 We say $a_n = o(b_n)$ if $\lim\sup_n \frac{a_n}{b_n} = 0$, $a_n = O(b_n)$ if $\lim\sup_n \frac{a_n}{b_n} = K < \infty $,  $\Theta(a_n) = b_n$ if $a_n = O(b_n)$ and $b_n = O(a_n)$.   $\delta_x$ indicates a  unit point mass at $x \in \cX$.
We write for $x, y \in \mathbb{R}$, $x \wedge y \defeq \min(x,y)$, $x \vee y \defeq \max(x,y)$. We denote the space of probability distributions over the finite set $\cX$ by $\Delta_{\cX}$. We reserve the symbols $k_x = \abs{\cX}$ and $k_y = \abs{\cY}$. We use the upper case of a letter to denote a random variable and the lower case to indicate the realization of that random variable. 
 Let $\loss: \Delta_{\cX} \times \Delta_{\cX} \rightarrow  \mathbb{R}$ be a loss function. 

In the minimax setting \cite{Wald49}, we assume that nature adversarially chooses a distribution $p_X$; $n$ samples $u^n = \{x_i\}^n_{i=1}$ are drawn i.i.d. from this distribution; and the goal is to design an estimator $\hat{q}: \cX^n \to \Delta_\cX$ based on the samples $u^n$ to minimize the expected loss. We denote the associated risk by $r^{\loss}_n$, and the problem is formulated as:  
\begin{equation}
    r^{\loss}_n \defeq \min_{\hat{q}_X} \max_{p_X \in \Delta_{\cX}} \Ed{U^n}{\loss(p_X, \hat{q}_X(U^n))} \label{eq:minimax_general}
\end{equation}
In this work we will consider two different forms of loss function:
\subsection{$l^p_p$ losses for $1 \le p < 2$}
In this case $\cL$ is chosen to be:
\begin{equation}
\loss(p_X,q_X) = \|p_X - q_X\|^p_p \defeq \sum_{x \in \cX} \Big( \abs{ p_X(x) - q_X(x) } \Big)^p
\end{equation}
In particular the risk  becomes for $l^p_p$ losses:
\begin{equation}
r^p_n \defeq  \min_{\hat{q}_X } \max_{p_X} \Ed{U^n}{ \|p_X - \hat{q}_X(U^n)\|^p_p} \label{eq:rmn2}
\end{equation}
The choices of $l^2_2$ and $l_1$ are has been considered in the prior work \cite{trybula1958some, wilczynski1985minimax, HanJW14}. For general $p$ as we consider in this work the minimax risk scales as $r^p_n = \Theta(n^{-\frac{p}{2}})$ \cite{onsemsup}. Similarly we will adopt the convention in \cite{onsemsup} and denote the constant of this rate by $r^n_p =  C_p (1 + o(1)) n^{-\frac{p}{2}}) $. 

One of the main goals of this paper is to study the extension of \eqref{eq:rmn2} into a multivariate semi-supervised setting. In this new setting, the nature chooses a model $p_{XY}$ of a pair of jointly distributed random variables $X,Y$. The estimator has access to two datasets: a collection of complete $l^m \defeq \{(x'_i, y'_i)\}^{m}_{i=1}$ and a collection of incomplete samples  $u^n \defeq \{x_i\}^n_{i=1}$ generated i.i.d from $p_{XY}$ and $p_X$. Our goal is to design an estimator with minimal expected risk.  We formulate this new problem as: 
\begin{IEEEeqnarray}{rCl} \label{eq:Rmn}
R^p_{m,n} \defeq \min_{\hat{q}_{XY}}\max_{p_{XY} } \Ed{U^n, L^m}{ \|p_{XY} - \hat{q}_{XY}(U^n, L^m)  \|^p_p}  \IEEEeqnarraynumspace
\end{IEEEeqnarray}
Throughout our study of \eqref{eq:Rmn}, we will consider the following auxiliary problems:
\begin{gather*}
    R^p_m \defeq \min_{\hat{q}_{Y\mid X }}\max_{p_{XY}} \Ed{L^m}{\|p_{XY} - p_X \hat{q}_{Y \mid X }(L^m)\|^p_p}    \\
\Bar{R}^p_m \defeq \min_{\hat{q}_{Y \mid X}} \max_{p_{X}} \Ed{L^m_{X}}{ \max_{p_{Y\mid X}} \Ed{L^m_Y}{\|p_{XY} - p_X \hat{q}_{Y \mid X}(L^m)\|^p_p}} 
\end{gather*}

\subsection{$f$-divergences}
$f$-divergences represent a crucial family of statistical loss functions. In the context of discrete distributions, f-divergences assume the form:
\begin{equation}
\loss(p_X,q_X) = D_f(p_X \| q_X) \defeq  \sum_{x \in \cX} q_X(x) f\left(\frac{p_X(x)}{q_X(x)}\right)
\end{equation}
for a function  strictly convex $f: \cX \rightarrow \reals$ with $f(1)=0$. The minmax estimation problem under consideration for this loss is speacialized to:
\begin{equation}
r^f_n \defeq \min_{\hat{q}_X } \max_{p_X \in \Delta_{\delta}} \Ed{U^n}{ D_f( p_X \|  \hat{q}_X(U^n) ) } \label{eq:rf}
\end{equation}
where $\Delta_{\delta} = \{ p_X \in \cP(\cX) : p_X(x) \ge \delta, \quad \forall x \in \cX \}$ and we further assume that $f$ is thrice differentiable and is subexponential, i.e. $\limsup_{x \rightarrow \infty} \frac{\abs{f(x)}}{e^{c x }} = 0$ for all $c \in \reals$. Under these constraint minmax risk is $r^f = C_f n^{-1} + o(n^{-1})$ $C_f \defeq  f''(1) (\cardx - 1) /2$ \cite{pmlr-v40-Kamath15}. In addition we impose on $f$ that $\lim_{x \rightarrow 0^+} \abs{f(x)} < \infty$ and $\lim_{x \rightarrow 0^+} x \abs{ f'(x)} < \infty$.  In the sequel, we will hide the underlying set $\Delta_\delta$.  We note that these mild assumptions are satisfied by common f-divergences including KL,$\chi^2$,Squared Hellinger and Le Cam divergences.

Similarly we define the main problem of interest for $f$-divergences:
\begin{IEEEeqnarray}{rCl} \label{eq:Rfmn}
R^f_{m,n} \defeq \min_{\hat{q}_{XY}}\max_{p_{XY} } \Ed{U^n, L^m}{ D_f(p_{XY} \| \hat{q}_{XY}(U^n, L^m) ) }  \IEEEnonumber\IEEEeqnarraynumspace
\end{IEEEeqnarray}
and auxilary problems:
\begin{gather*}
    R^f_m \defeq \min_{\hat{q}_{Y\mid X }}\max_{p_{XY}} \Ed{L^m}{ D_f( p_{Y\mid X} \|  \hat{q}_{Y \mid X }(L^m) \mid p_X ) }    \\
\Bar{R}^f_m \defeq \min_{\hat{q}_{Y \mid X}} \max_{p_{X}} \mathbb{E}_{L^m_{X}}{ \max_{p_{Y\mid X}} \Ed{L^m_Y}{ D_f( p_{Y\mid X} \| \hat{q}_{Y \mid X} (L^m) \mid p_X) }} 
\end{gather*}
where we introduce the \textit{conditional $f$-divergence}:
\begin{equation*}
D_f(p_{Y\mid X} \| q_{Y\mid X} | p_X) = \sum_{x \in \cX } p_X(x) D_f( p_{Y\mid X =x }\| q_{Y\mid X = x} ) 
\end{equation*}

\section{Results}
\subsection{$l^p_p$ loss functions}
\begin{thm}[Theorem 1 of \cite{onsemsup}]  \label{thm:1} 
Let $\hat{q}^*_{n}$ be a minimax optimal estimator for $r^p_n$. Then the conditional composition $\hat{q}^{*,m}_{Y\mid X}$ based on $ \hat{q}^{*}_n$ is minimax optimal for $\bar{R}^p_m$:
\begin{IEEEeqnarray}{lCr}
\max_{p_{X}} \Ed{L^m}{\max_{p_{Y|X}}\|p_{XY} - p_X \hat{q}^{*,m}_{Y \mid X} \|^p_p} = \bar{R}^p_m 
\end{IEEEeqnarray}
\end{thm}

\begin{thm}\label{thm:lp:2}
Let $1 \le p\le 2$ and  $\hat{q}^*_n$ be a first order minimax optimal estimator for $r^p_n$. Then the conditional composition $\hat{q}^{*,m}_{Y\mid X}$ based on $\hat{q}^{*}_n$ is first order minimax optimal for $R^p_m$:
\begin{IEEEeqnarray}{lCr}
 \max_{p_{XY}} \Ed{L^m}{ \|p_{XY} - p_X \hat{q}^{*,m}_{Y\mid X} \|^p_p }  = R^p_m + o\left( R^p_m \right)
\end{IEEEeqnarray}
\end{thm}

\begin{thm}[Theorem 3 of \cite{onsemsup}]  \label{thm:3} Let $m = o(n)$:
\begin{IEEEeqnarray}{rCl}
\abs{ R^p_{m,n} - R^p_m   } \le  O\left( {m^{-\frac{p -1 }{2}} (n)^{-1/2} } \right)
\end{IEEEeqnarray}
\end{thm}

\begin{thm} \label{thm:lp_main}
Let $1 \le p \le 2$ and $\hat{q}^*_n$ be a first-order optimal estimator for $r^p_n$ with minimax risk: 
\begin{align*}
    \max_{p_{X}} \Ed{}{\|p_X -\hat{q}^*_n\|^p_p} = \min_{\hat{q}} \max_{p_{X}} \Ed{}{\|p_X -\hat{q}_n\|^p_p} \simeq\frac{C_p}{n^{\frac{p}{2}}} 
\end{align*}
Then the joint composition $\hat{q}^{*, m,n}_{XY}$ based on $\hat{q}_n$ is first order minimax optimal for $R^p_{m,n}$ in the regime $m=o(n).$ \\ Furthermore,
\begin{align*}
 R^p_{n,m} =  \left( k_x \right)^{1 - \frac{p}{2}} \frac{C_p}{m^{\frac{p}{2}}} + \littleOh{\frac{1}{m^{\frac{p}{2}}}}
\end{align*}
\end{thm}
\subsection{f-divergences}
\begin{thm}\label{thm:f:1} 
Let $\hat{q}^*_{n}$ be a minimax optimal estimator for $r^f_n$. Then the conditional composition $\hat{q}^{*,m}_{Y\mid X}$ based on $ \hat{q}^{*}_n$ is minimax optimal for $\bar{R}^f_m$:
\begin{IEEEeqnarray}{lCr}
\max_{p_{X}} \Ed{L^m}{\max_{p_{Y|X}} D_f( \check{q}^{*,m}_{Y\mid X} \| \check{q}^*_n \mid p_X) } = \bar{R}^f_m 
\end{IEEEeqnarray}
\end{thm}
\begin{corol} \label{thm:f:2}
Let $\hat{q}^*_n$ be a first order minimax optimal estimator for $r^f_n$. Then the conditional composition $\hat{q}^{*,m}_{Y\mid X}$ based on $\hat{q}^{*}_n$ is first order minimax optimal for $R^p_m$:
\begin{IEEEeqnarray}{lCr}
 \max_{p_{XY}} \Ed{L^m}{ D_f( p_{Y\mid X} \| \hat{q}^{*,m}_{Y\mid X}  \mid p_X) }  = R^f_m + o\left( R^f_m \right)
\end{IEEEeqnarray}
\end{corol}
\begin{thm} \label{thm:f:3}
Let $m=o(n)$:
$$
\abs{R^f_{m,n} - R^f_m} = o(R^f_m)
$$
\end{thm}
\begin{thm} \label{thm:f:main}
Let $\hat{q}^*_n$ be a first-order optimal estimator for $r^p_n$ with minimax risk: 
\begin{align*}
    \max_{p_{X}} \Ed{}{ D_f( p_X \| \hat{q}^*_n ) } = \min_{\hat{q}} \max_{p_{X}} \Ed{}{ D_f(p_X \| \hat{q}_n) } \simeq C_f n^{-1} 
\end{align*}
Then the joint composition $\hat{q}^{*, m,n}_{XY}$ based on $\hat{q}_n$ is first order minimax optimal for $R^p_{m,n}$ in the regime $m=o(n).$ \\ Furthermore,
\begin{align*}
 R^f_{n,m} = \cardx \frac{C_f}{m} + o\left(\frac{1}{{m}}\right) 
\end{align*}
\end{thm}

\section{Proofs for Theorems}
\subsection{$l^p_p$ loss functions}
\begin{IEEEproof}[Proof for \thmref{thm:lp:2}]
By (12) in the proof for \thmref{thm:1},
\begin{IEEEeqnarray}{rCl}
\bar{R}^p_m &=& \max_{p_X} \sum_{x \in \cX} \sum^{m}_{i=0} {m \choose i} \left( p_X(x) \right)^{i+p} \left( 1 - p_X(x) \right)^{m-i} r^p_i \IEEEnonumber \\ 
\label{eq:lem:bar_red_lp:1} \\
&=& \max_{p_X} \sum_{x \in \cX} C_p \left( \frac{p_X(x)}{m} \right)^{\frac{p}{2}}   + \littleOh{\frac{1}{m^{\frac{p}{2}}}} \label{eq:lem:bar_red_lp:2} \\ 
&=& \sum_{x \in \cX} C_p \left( \frac{\frac{1}{\abs{\cX}}}{m} \right)^{\frac{p}{2}} +  \littleOh{\frac{1}{m^{\frac{p}{2}}}}\label{eq:lem:bar_red_lp:3}  \\
&=& {\abs{\cX}}^{1 - \frac{p}{2}} \frac{C_p}{m^{\frac{p}{2}}} + \littleOh{\frac{1}{m^{\frac{p}{2}}}} \label{eq:lem:bar_red_lp:4} \nonumber 
\end{IEEEeqnarray}
\eqref{eq:lem:bar_red_lp:2} follows from \lemref{lem:convergence_final}.  In \eqref{eq:lem:bar_red_lp:1}, the objective function is symmetric in variables $\{p_X(x) : x \in \cX\}$ and for $p \le 2$ each summand is a concave function. Therefore the optimal solution should be the uniform distribution, i.e. $p^*(x) = \frac{1}{\abs{\cX}}$  which leads to  \eqref{eq:lem:bar_red_lp:3}.

On the other hand, 
\begin{IEEEeqnarray}{rCl}
R^p_m &\defeq& \min_{\hat{q}_{Y\mid X}} \max_{p_{XY}} \Ed{L^m}{\|p_{XY} - p_X \hat{q}_{Y\mid X}(L^m)\|^p_p} \label{eq:lem:bar_red_lp:2:1}\\ 
&\ge&  \min_{\hat{q}_{Y\mid X}}\max_{p_{Y\mid X}} \Ed{L^m}{\| U_{\cX} p_{Y\mid X} - U_{\cX} \hat{q}_{Y\mid X}(L^m)\|^p_p} \label{eq:lem:bar_red_lp:2:2}  \\ 
&=&  \min_{\hat{q}_{Y\mid X}}\max_{p_{Y\mid X}} \Ed{L^m}{\sum_{x \in \cX} \frac{1}{\cardx^p} \|p_{Y \mid X=x} - \hat{q}_{Y\mid X=x}\|^p_p }  \label{eq:lem:bar_red_lp:2:3}\\
&=& \min_{\hat{q}_{Y\mid X}}  \sum_{x \in \cX} \frac{1}{\cardx^p} \max_{p_{Y\mid X=x}} \Ed{L^m}{\|p_{Y \mid X=x} - \hat{q}_{Y\mid X=x}\|^p_p } \IEEEeqnarraynumspace   \label{eq:lem:bar_red_lp:2:3} \\
&=& \sum_{x \in \cX}\frac{1}{\cardx^p} \min_{\hat{q}_{Y\mid X=x}}  \max_{p_{Y\mid X=x}}  \Ed{L^m}{\|p_{Y \mid X=x} - \hat{q}_{Y\mid X=x}\|^p_p }   \IEEEnonumber \\
\label{eq:lem:bar_red_lp:2:4} \\
&=& \sum_{x \in \cX} \frac{1}{\cardx^p} (1 + o(1)) \frac{C_p}{ ( \frac{m}{\cardx} )^{\frac{p}{2}}}  \label{eq:lem:bar_red_lp:2:5}\\
&=& (1 + o(1) ) \cardx^{1 -\frac{p}{2}} \frac{C_p}{m^{\frac{p}{2}}} \label{eq:lem:bar_red_lp:2:6}
\end{IEEEeqnarray}
In \eqref{eq:lem:bar_red_lp:2:1} we substitute uniform distribution as $p_X$ in the inner maximization problem to arrive at the lower bound in \eqref{eq:lem:bar_red_lp:2:2}. In \eqref{eq:lem:bar_red_lp:2:3} we expand the $l^p_p$ norm. In $\eqref{eq:lem:bar_red_lp:2:4}$ we note that optmization variables $\{p_{Y\mid X=x} : x \in \cX \}$ and $\{\check{q}_{Y\mid X=x} : x \in \cX \}$ are independent in the order. \eqref{eq:lem:bar_red_lp:2:5} follows from \lemref{lem:hoeffding_lemma}. 
Finally, lemma follows by the relation $R^p_m \le \bar{R}^p_m $ along with \eqref{eq:lem:bar_red_lp:4} and \eqref{eq:lem:bar_red_lp:2:6}
\end{IEEEproof}

\begin{IEEEproof}[Proof for \thmref{thm:lp_main}]
By Theorem-3 of \cite{onsemsup} $\abs{ R^p_{m,n } - R^p_m } \le O(m^{-\frac{p-1}{2}} n^{-\frac{1}{2}})$. Furthermore, by Theorem-1 the composition estimator is minimax optimal estimator for $\bar{R}^p_m$. However by \thmref{thm:lp:2} $\bar{R}^p_m - R^p_m \le o\left(R^p_m\right)$ for $1 \le p < 2$ as well. Therefore, the composition estimator is first order minimax for $R^p_{m,n}$ for $1\le p \le 2$. 
\end{IEEEproof}

\subsection{$f$-divergences}
\begin{IEEEproof}[Proof for \thmref{thm:f:1}]
Let us define:
\begin{IEEEeqnarray}{C}
    h(p_X, \hat{q}_{Y\mid X}) \defeq \Ed{L^m_X}{ \max_{p_{Y\mid X}} \Ed{L^m_Y}{ D_f( p_{Y\mid X} \| \hat{q}_{Y\mid X}  \mid p_X ) } } \IEEEnonumber \\
    \IEEEeqnarraynumspace \label{eq:defn_eq_fpx}
\end{IEEEeqnarray}
and let $\hat{q}^{**}_{Y\mid X} : (\cX \times \cY)^m \rightarrow (\Delta_{\cY})^{\abs{\cX}}$ be an  estimator for the conditional distribution $p_{Y\mid X}$. Because of the inequalities $\max_{p_X} h(p_X, \hat{q}^{**}_{Y \mid X}) \ge \bar{R}^p_m \defeq \min_{\hat{q}} \max_{p_X} \ h(p_X, \hat{q}) \ge  \max_{p_X} \min_{\hat{q} } h(p_X, \hat{q})$ it is sufficient  to show that for all $p_X$, $\min_{\hat{q}_{Y\mid X}} h(p_X, \hat{q}_{Y\mid X} ) = h(p_X,  \hat{q}^{**}_{Y\mid X} )$. Now let us show the latter is satisfied by  the composition estimator $\hat{q}^{*,m}_{Y \mid X}$. Set a $p_X$ and we define $p_{i,x} \defeq \prob{   T_x(L_X) = i}$. Then the left-hand side becomes:
\begin{IEEEeqnarray*}{rCl}
\min_{\hat{q}_{Y\mid X}} h(p_X, \hat{q}_{Y\mid X} )  &=& \min_{\hat{q}_{Y\mid X}} \Ed{L^m_Y}{ D_f( p_{Y\mid X} \| \hat{q}_{Y\mid X}  \mid p_X ) }\\
\IEEEeqnarraymulticol{3}{l}{= \min_{\hat{q}_{Y\mid X}}  \mathbb{E}_{L^m_X} \Big[ {\max_{p_{Y\mid X}} \mathbb{E}_{L^m_Y} \Big[ {\sum_{x \in \cX }  p_x D_f( p_{Y | X = x } \|  \hat{q}_{Y\mid X=x } ) }\Big] } \Big]  }  \\
  \IEEEeqnarraymulticol{3}{l}{ =\sum_{x \in \cX } p_x
 \min_{\hat{q}_{Y|X=x}}  \mathbb{E}_{L^m_X} \Big[ {\max_{p_{Y|X = x}} \Ed{L^m_Y}{ D_f( p_{Y | X = x } \|   \hat{q}_{Y | X=x } ) } }  \Big] }  \IEEEeqnarraynumspace\\ 
 \IEEEeqnarraymulticol{3}{l}{\IEEEyesnumber \label{eq:risk_substitute_pre} }  \\
 \IEEEeqnarraymulticol{3}{l}{ =\sum_{x \in \cX }  p_X(x)  \Big( \sum^m_{i=0} p_{i,x}    r^f_i  \IEEEyesnumber \label{eq:risk_substitute} \Big)}
\end{IEEEeqnarray*}
In \eqref{eq:risk_substitute_pre} we note that  the optimization variables  are independent. Steps leading \eqref{eq:risk_substitute_pre} to \eqref{eq:risk_substitute} are given below and these steps demonstrate that $\min_{\hat{q}_{Y\mid X}} h(p_X, \hat{q}_{Y\mid X} ) = h(p_X,  \hat{q}^{**}_{Y\mid X} )$ holds for the estimator $\hat{q}^{*,m}_{Y \mid X}$:
\begin{IEEEeqnarray*}{rCl}
    \IEEEeqnarraymulticol{3}{l}{  \min_{\hat{q}_{Y\mid X=x}}  \Ed{L^m_X}{\max_{p_{Y\mid X = x}} \Ed{L^m_Y}{ D_f( p_{Y \mid X = x } \|  \hat{q}_{Y\mid X=x } ) } } } \\
     &=&    \min_{\hat{q}_{Y\mid X=x}}  \sum^m_{i=0}   p_{i,x}  \max_{p_{Y\mid X = x}} \Ed{L^i_Y}{D_f(p_{Y\mid X =x } \| \hat{q}^i_{Y\mid X =x } )}   \\ 
     \IEEEeqnarraymulticol{3}{l}{
     = \sum^m_{i=0}  p_{i,x}  \min_{\hat{q}^i_{Y\mid X=x}} \max_{p_{Y\mid X = x}}  \Ed{L^i_Y}{ D_f( p_{Y\mid X =x } \| \hat{q}^i_{Y\mid X =x } )}  } \IEEEeqnarraynumspace \IEEEyesnumber \label{eq:risk:another} \\
    \IEEEeqnarraymulticol{3}{l}{
     = \sum^m_{i=0}  p_{i,x}   \max_{p_{Y\mid X = x}}  \Ed{L^i_Y}{ D_f( p_{Y\mid X =x } \| \hat{q}^{*,i}_{Y\mid X =x } )}  } \IEEEeqnarraynumspace \IEEEyesnumber \label{eq:risk:another_post}
\end{IEEEeqnarray*}
To arrive in \eqref{eq:risk:another_post}, we notice in \eqref{eq:risk:another} that  $\min_{\hat{q}^i_{Y\mid X=x}} \max_{p_{Y\mid X = x}}  \Ed{L^i_Y}{D_f(p_{Y\mid X =x }\| \hat{q}_{Y\mid X =x } )}$ is the problem $r^f_i$ and hence is achieved by $\hat{q}^{*,i}_{Y \mid X}$. 
\end{IEEEproof}
\begin{IEEEproof}[Proof for \corolref{thm:f:2}]
We proceed similar to the proof for \thmref{thm:lp:2}.
We have by \eqref{eq:risk_substitute}:
\begin{align}
    \bar{R}^f_m =& \sum_x p_x \sum^m_{i=0} {m \choose i } p^i_x (1-p_x)^{m - i} r^f_i \nonumber \\ 
    =&  \sum_x \frac{1}{p_x} \sum^m_{i=0} {m \choose i }  p^{2+i}_x (1-p_x)^{m - i} r^f_i \label{eq:use_the_prev_1} \\ 
    \simeq& \sum_x \frac{1}{p_x} C_f \frac{p_x}{m}=  \cardx \frac{C_f}{m} \label{eq:use_the_prev_2}
\end{align}
In \eqref{eq:use_the_prev_2} we evaluate the sum \eqref{eq:use_the_prev_1} via \lemref{lem:convergence_final}.Similar to the proof for \thmref{thm:lp:2}, we have the lower bound for $R^f_m$:
\begin{align}
  R^f_m =& \min_{\cq_{Y\mid X}} \max_{p_{XY}} \mathbb{E}_{L^m}[ \sum_x p_x D_f(p_{Y\mid X=x} \| \cq_{Y\mid X=x}) ] \nonumber \\  
  \ge& \min_{\cq_{Y\mid X}} \max_{p_{Y \mid X}} \mathbb{E}_{L^m}[ \sum_x \frac{1}{\cardx} D_f(p_{Y\mid X=x} \| \cq_{Y\mid X=x}) ] \label{eq:unif_subs} \\ 
  =& \sum_{x \in \cX}  \frac{1}{\cardx} \min_{\cq_{Y\mid X=x}} \max_{p_{Y \mid X=x}} \mathbb{E}_{L^m}[  D_f(p_{Y\mid X=x} \| \cq_{Y\mid X=x}) ]\nonumber \\ 
  =& \sum_{x \in \cX} \frac{1}{\cardx} \frac{C_f}{ \frac{m}{\cardx} } (1 + o(1)) \label{eq:unif_subs:3} \\ 
  =& \cardx \frac{C_f}{m} (1 + o(1))
\end{align}
in \eqref{eq:unif_subs} we substitute $p_X = \cU_{\cX}$. In \eqref{eq:unif_subs:3} we use \lemref{lem:hoeffding_lemma} with $p=2$ as $r^f_i = \Theta(i^{-1})$. We  omit a version of this lemma for f-divergences as the same proof applies.  Finally we note that $R^f_m \le \bar{R}^f_m$.
\end{IEEEproof}

\begin{IEEEproof}[Proof for \thmref{thm:f:3}]
We have the lower bound: 
$$
R^f_{n+1, m} \le R^f_{n,m} \implies R^f_m = \lim_{n \rightarrow \infty} R^f_{m,n} \le R^f_{m,n} 
$$
Now let us prove the upperbound: 
\begin{IEEEeqnarray*}{rCl}
R^f_{n,m} &=& \min_{\hat{q}_{XY}} \max_{p_{XY}} \Ed{U_n,L_m}{D_f(p_{XY} \| \cq_{XY}(U_n, L_m) )}  \\
&=& \min_{\hat{q}_{XY}} \max_{p_{XY}} \Ed{U_n,L_m}{ \sum_{x,y} \hat{q}_{XY}(x,y) f\left( \frac{p_{XY}(x,y)}{\cq_{XY}(x,y)} \right) } \\
\IEEEeqnarraymulticol{3}{l}{= \min_{\hat{q}_{XY}} \max_{p_{XY}} \mathbb{E}_{U_n, L_m} \Big[ \sum_{x,y} p_X(x) \hat{q}_{Y\mid X}(y \mid x) f\left( \frac{p_{Y\mid X}(y \mid x)}{\hat{q}_{Y\mid X}(y \mid x)} \right)  }  \\ 
&\quad&+  (p_X(x) - \hat{q}_X(x) ) \Gamma^{U_n, L_m}_{x,y}  \Big]
\label{eq:taylor_expansion}  \IEEEyesnumber\\
\IEEEeqnarraymulticol{3}{l}{\le \min_{\hat{q}_{XY}} \max_{p_{XY}} \mathbb{E}_{U_n, L_m} \Big[ \sum_{x,y} p_X(x) \hat{q}_{Y\mid X}(y \mid x) f\left( \frac{p_{Y\mid X}(y \mid x)}{\hat{q}_{Y\mid X}(y \mid x)} \right) \Big] } \\
&+& \min_{\hat{q}_{XY}} \max_{p_{XY}} \mathbb{E}_{U_n, L_m} \Bigg[ \sum_{x,y} (p_X(x) - \hat{q}_X(x) )  \hat{\Gamma}_{x,y} \Bigg] \\
\IEEEeqnarraymulticol{3}{l}{= R^f_m + \min_{\hat{q}_{XY}} \max_{p_{XY}} \mathbb{E}_{U_n, L_m} \bigg[ \sum_{x,y} (p_X(x) - \hat{q}_X(x) )  \Gamma^{U_n, L_m}_{x,y} \bigg] } \\
\IEEEeqnarraymulticol{3}{l}{ \le R^f_m + \min_{\hat{q}_{XY}} \max_{p_{XY}} \mathbb{E}\Big[{ \Big( { {\sum_{x} (p_x - \cq_x)^2 } } \Big)^{\frac{1}{2}} \cdot \Big({ {\sum_{x} (\sum_y \hat{\Gamma}_{x,y})^2  } } \Big)^{\frac{1}{2}} } \Big]    } \\
\end{IEEEeqnarray*}
\begin{IEEEeqnarray*}{rCl}
\IEEEeqnarraymulticol{3}{l}{\le R^f_m + \min_{\hat{q}_{XY}} \max_{p_{XY}} \mathbb{E} \Big[{ \sum_{x} (p_x - \cq_x)^2  }\Big] \mathbb{E}\Big[{  { {\sum_{x} (\sum_y \hat{\Gamma}_{x,y})^2  } } } \Big] }\\
\IEEEeqnarraymulticol{3}{l}{\le R^f_m + }\\ 
\IEEEeqnarraymulticol{3}{l}{  \min_{\hat{q}_{XY}} \max_{p_{XY}} \mathbb{E}\Big[{ \sum_{x} (p_x - \cq_x)^2  } \Big] \cdot \min_{\hat{q}_{XY}} \max_{p_{XY}} \mathbb{E}\Big[{  { {\sum_{x} (\sum_y \hat{\Gamma}_{x,y})^2  } } }\Big] }    \\
 \IEEEeqnarraymulticol{3}{l}{= R^f_m +  r^2_n  \cdot \min_{\hat{q}_{XY}} \max_{p_{XY}} \mathbb{E}\Big[{  { {\sum_{x} (\sum_y \hat{\Gamma}_{x,y})^2  } } } \Big]  } \\
&\le& R^f_m +  r^2_n  \cdot  \max_{p_{XY}} \mathbb{E}\Big[{  { {\sum_{x} (\sum_y \Gamma_{x,y})^2  } } }\Big] \label{eq:beta_subs} \IEEEyesnumber 
\end{IEEEeqnarray*}
In \eqref{eq:taylor_expansion}, we apply mean value theorem pointwise to $p_X \rightarrow p_{XY} f\left( \frac{p_{XY}(x,y)}{\hat{q}_{XY}(x,y)}\right)$ where the remainder is evaluated at a point $\xi^{u_n, l_m}_{x,y} \in (p_X(x), \cq_X(x))$. For the compactness of the notation we abbreviate gradient terms with:
$$
\Gamma^{U_n, L_m}_{x,y} \defeq \nabla_{q} \left(  q(x) \cq_{Y\mid X}(y\mid x) f\left(\frac{p_{XY}(y \mid x) }{ \xi^{U_n,L_m}_{y,x} \cq_{Y\mid X}(y\mid x)} \right)  \right) 
$$
and in \eqref{eq:beta_subs} we chose the estimator $\hat{q}_{XY}(x,y) = \frac{1}{\cardx \cdot \cardy}$ for the outer minimization problem which we denote by $\Gamma_{x,y}$.  Observe that the term $\sqrt{r^2_n} = \Theta(n^{-\frac{1}{2}})$.  Therefore we are done if its multiplier $\min_{\hat{q}_{XY}} \max_{p_{XY}} \E{  { {\sum_{x} (\sum_y \hat{\Gamma}_{x,y})^2  } } } $ is bounded. In the sequel we show this via the pointwise boundedness of $\hat{\Gamma}_{x,y}$ which is the case under our assumptions on $f$. First notice that by the convexity of the perspective function the map $p_X \rightarrow p_{XY} f\left( \frac{p_{XY}(x,y)}{\hat{q}_{XY}(x,y)}\right)$ is convex and thereby its gradients are monotone. Hence it can be bounded by its values at its end points where either $\xi^{U_n, L_m}_{y,x} = p_X(x) $ or $\xi^{U_n, L_m}_{y,x} =\cq_X(x) = \frac{1}{\cardx}$, therefore we analyze these two cases seperately:
\paragraph{$\xi^{U_n, L_m}_{y,x} = p_X(x)$}
\begin{IEEEeqnarray*}{rCl}
\IEEEeqnarraymulticol{3}{l}{ \abs{\hat{\Gamma}_{x,y}} = \Big| \frac{1}{\cardy} f\Big( \frac{p_{Y | X}(y |  x)}{\frac{1}{\cardy}} \Big) - \frac{p_{Y| X}(y | x)}{\frac{1}{\cardy}} f'\Big( \frac{p_{Y| X}(y| x) }{\frac{1}{\cardy}} \Big) \Big|
} \\
\IEEEeqnarraymulticol{3}{l}{
\le \Big| \frac{1}{\cardy} f\Big( \frac{p_{Y | X}(y | x)}{\frac{1}{\cardy}} \Big) \Big| + \Big| \frac{p_{Y | X}(y | x)}{\frac{1}{\cardy}} f'\Big( \frac{p_{Y | X}(y  x) }{\frac{1}{\cardy}} \Big) \Big|
} \\\IEEEyesnumber \label{eq:unbounded_case:1}
\end{IEEEeqnarray*}
In \eqref{eq:unbounded_case:1} can be unbounded only when $p_{Y \mid X} \rightarrow 0$, which creates the cases $f(0^+)$ and $0^+ f'(0^+)$ which by assumption are bounded.
\paragraph{$\xi^{U_n, L_m}_{y,x} =\cq_X(x) = \frac{1}{\cardx}$}
Similar to (a), after substituting the estimator we arrive at:
\begin{IEEEeqnarray*}{rCl}
\IEEEeqnarraymulticol{3}{l}{  \abs{ \hat{\Gamma}_{x,y} } =
\bigg| \frac{1}{\cardy} f\bigg( \frac{p_{XY}(x,y)}{ \frac{1}{\cardx \cardy} } \bigg)  
- \frac{p_{XY}(x,y)}{ \frac{1}{\cardx}} f'\bigg( \frac{p_{XY}(x,y)}{\frac{1}{\cardx \cardy} } \bigg)  \bigg| } \\
\IEEEeqnarraymulticol{3}{l}{ \le 
\bigg| \frac{1}{\cardy} f\bigg( \frac{p_{XY}(x,y)}{ \frac{1}{\cardx \cardy} } \bigg) \bigg|
+ \bigg| \frac{p_{XY}(x,y)}{ \frac{1}{\cardx}} f'\bigg( \frac{p_{XY}(x,y)}{\frac{1}{\cardx \cardy} } \bigg)  \bigg| } 
\\
\IEEEyesnumber
\label{eq:unbounded_case:2}
\end{IEEEeqnarray*}
\eqref{eq:unbounded_case:2} can be unbounded only when $p_{XY}(x,y) \rightarrow 0$, which again creates the cases $f(0^+)$ and $0^+ f'(0^+)$ which by assumption are bounded.
\end{IEEEproof}
\begin{IEEEproof}[Proof for \thmref{thm:f:main}]
By \thmref{thm:f:3} $\abs{ R^f_{m,n } - R^f_m } \le o(R^f_m)$. Furthermore, by \thmref{thm:f:1} the composition estimator is minimax optimal estimator for $\bar{R}^f_m$. On the other hand by \thmref{thm:f:2} $\bar{R}^f_m - R^f_m \le o\left(R^p_m\right)$. Therefore, the composition estimator is first order minimax for $R^f_{m,n}$. 
\end{IEEEproof}
\section{Supplementary Results}
\begin{lemma} \label{lem:hoeffding_lemma}For all $x \in \cX$: 
$$
\min_{\hat{q} } \max_{p_{Y\mid X=x}} \Ed{ L^m \sim \cU_{\cX} p_{Y\mid X} }{\|p_{Y \mid X= x} - \hat{q}\|^p_p} \simeq C_p \cdot \Big( \frac{\cardx}{m}\Big)^{\frac{p}{2}}
$$
\end{lemma}
\begin{IEEEproof}
Let us fix  constants $C > 0$ and $\frac{1}{2} < \alpha < 1$.Let us denote the event $\cA \defeq \{ \frac{m}{\cardx} - C \cdot m^{\alpha} \le T_x \le \frac{m}{\cardx} + C \cdot m^{\alpha}\}$ and we let $s_m = m \vee ( C \cdot m^{\alpha} + m/\cardx )$ Observe that:
\begin{IEEEeqnarray}{rCl}
\IEEEeqnarraymulticol{3}{l}{ \Ed{L^m}{\|p_{Y \mid X=x} - \hat{q}\|^p_p } \ge\prob{ \cA } \cdot \E{ \|p_{Y \mid X=x} - \hat{q}\|^p_p  \mid  \cA } }  \IEEEnonumber \\
\IEEEeqnarraymulticol{3}{l}{\ge\prob{ \cA } \E{ \|p_{Y \mid X=x}- \hat{q}_{Y\mid X=x}\|^p_p  \mid  T_x =  s_m } \label{eq:lem:hoeffding_lemma:1} } \IEEEeqnarraynumspace   \label{eq:lem:hoeffding_lemma:1}\\
\IEEEeqnarraymulticol{3}{l}{ \ge (1 - 2 e^{-C^2 m^{2 \alpha - 1 }} )\cdot \E{ \|p_{Y \mid X=x}- \hat{q}_{Y\mid X=x}\|^p_p  \mid  T_x =  s_m }  } \IEEEeqnarraynumspace \label{eq:lem:hoeffding_lemma:2}
\end{IEEEeqnarray}
In \eqref{eq:lem:hoeffding_lemma:1} we lower bound the conditional expectation by noting that the risk is smaller if more samples are from the category $x$. In \eqref{eq:lem:hoeffding_lemma:2} we use bounded differences inequality to bound the probability of the event $\prob{\cA}$. Finally, taking the $\min\max$ of both sides we establish that:
\begin{gather*}
\min_{\hat{q} } \max_{p_{Y\mid X=x}} \Ed{ L^m \sim \cU_{\cX} p_{Y\mid X} }{\|p_{Y \mid X= x} - \hat{q}\|}  \\
\ge \left(1  - 2 e^{-C^2 m^{2 \alpha - 1 }}  \right) \cdot r^p_{ m \vee \frac{1}{\cardx} + C \cdot m^{\alpha} } \simeq C_p \cdot \Big( \frac{\cardx}{m}\Big)^{\frac{p}{2}} 
\end{gather*}
\end{IEEEproof}

We note that  \lemref{lem:convergence_final}, \lemref{lem:bernstein_trick} are identical to Lemma 1 and Lemma 2 of \cite{onsemsup} and except with adjustments tailored to make them work for the case $1 \le p \le 2$, which requires a more delicate analysis.
For  \lemref{lem:convergence_final} and \lemref{lem:bernstein_trick} we introduce:
\begin{IEEEeqnarray}{rCl}
H^n_p(x) &\defeq& \sum^n_{i=0} {n \choose i } r^p_i x^{i+p} (1- x)^{n- i} \\
    G^n_p(x) &\defeq& \sum^n_{i=0} \binom{n}{i} r^p_i x^{i} \left(\frac{i}{n} \right)^p (1-x)^{n-i} 
\end{IEEEeqnarray}

\begin{lemma} \label{lem:convergence_final} 
For $1 \le p < 2$: 
$$ 
 H^n_{p}(x) = C_p \ \left( \frac{x}{n} \right)^{\frac{p}{2}} +  o(n^{-\frac{p}{2}}) 
$$
\end{lemma}
\begin{IEEEproof}
Let $c$ be the constant indicated in the statement of \lemref{lem:bernstein_trick}. 
There are two cases to consider:
\begin{enumerate}
\item $x \ge c \frac{log^2(n)}{n}$:
\begin{IEEEeqnarray}{rCl}
H^n_{p}(x) &\simeq&  \sum^n_{i=0} \binom{n}{i}  \frac{C_p}{i^{\frac{p}{2}} + 1 }  x^{i+p} (1 -x )^{n -i}  \label{eq:lem1_1} \\
\IEEEeqnarraymulticol{3}{l}{= \frac{C_p}{n^{\frac{p}{2}}} \sum^n_{i=0} {n \choose i} \left(  \frac{i}{n} \right)^{\frac{p}{2}}  x^i (1- x)^{n -i}  + \bigOh{ \frac{ H^n_{p}(x)}{\sqrt{\log{n}}} } }  \IEEEnonumber \\  \label{eq:lem1_2} \\ 
\IEEEeqnarraymulticol{3}{l}{= \frac{C_p}{n^{\frac{p}{2}}} \left( x^{\frac{p}{2}} + O\left(n^{-\frac{p}{2}}\right) \right) + \bigOh{\frac{ H^n_{p}(x)}{\sqrt{\log{n}}} }   } \label{eq:lem1_3} \\
\IEEEeqnarraymulticol{3}{l}{= C_p \left( \frac{x}{n}\right)^\frac{p}{2}  + \littleOh{n^{-\frac{p}{2}}}}\label{eq:lem1_4}
\end{IEEEeqnarray}
\enlargethispage{-4cm}
\eqref{eq:lem1_1} holds since $r^p_n \simeq {C_p} {n^{-\frac{p}{2}}}$ whereas in \eqref{eq:lem1_2} we use \lemref{lem:bernstein_trick}. To obtain \eqref{eq:lem1_3}, we use \thmref{thm:approximation_bernstein} where  we let $f(x) = x^{\frac{p}{2}}$ and bound the error of $n$th order Bernstein polynomial approximation $B_n$ as: 
\begin{IEEEeqnarray*}{rCl}
   \IEEEeqnarraymulticol{3}{C}{  \abs{ B_n(x;f)- f(x)} = \frac{1}{n} \abs{ x (1-x) \frac{f''(x)}{2} } + o(n^{-1}) }  \\
    \IEEEeqnarraymulticol{3}{l}{  = \frac{1}{n} x (1-x) \frac{p}{4} (1 - \frac{p}{2}) x^{\frac{p}{2} - 2 } + o(n^{-1}) } \\
  \IEEEeqnarraymulticol{3}{l}{   \le \frac{1}{n^{\frac{p}{2}}}  \frac{p}{4} \Big(1 - \frac{p}{2} \Big) \left( {c \log^2(n) }\right)^{  \frac{p}{2} -1  } + o(n^{-1}) } \label{eq:step_step}  \IEEEyesnumber\\
\end{IEEEeqnarray*}
 where in \eqref{eq:step_step} we alleviate the bounds $1-x \le 1$ and $x \ge c \frac{log^2(n)}{n}$. Finally by \eqref{eq:lem1_3} we have $H^n_p(x) = O({n^{-\frac{p}{2}}})$ and we obtain \eqref{eq:lem1_4} by substituting it.
\item  $x \le c  \frac{log^2(n) }{n}$: 
\begin{IEEEeqnarray}{rCl}
H^n_{p}(x) &=&   \sum^n_{i=0} {n \choose i} r^p_i x^{i + p} (1- x )^{n -i} \\
&\le&  C_p \ x^p  \sum^n_{i=0} {n \choose i} \frac{1}{i^{\frac{p}{2}} + 1} x^i (1 - x)^{n - i} \\ 
&\le& C_p\ x^p = \bigOh{ \frac{\log^{2p} (n) }{n^p} }  
\end{IEEEeqnarray}
Similarly $ C_p \left(\frac{x}{n}\right)^{\frac{p}{2}} = O\left( \frac{\log^{2p}(n)}{n^p}  \right) $ when $x \le c \frac{  \log^2{n} }{n}$, hence $\abs{ C_p \left( \frac{x}{n} \right)^{\frac{p}{2}} - H^n_p(x) }  = \littleOh{n^{-\frac{p}{2}}} $. 
\end{enumerate}
\end{IEEEproof}

\begin{lemma}\label{lem:bernstein_trick}
There exists a $c > 0$ such that for $x \ge  c  \frac{log^2 n }{n}$: 
\begin{equation}
\abs{H^n_p(x) - G^n_p(x)} =  \bigOh{ \frac{  H^n_p(x) }{ \sqrt{\log{n}} } }
\end{equation}
\begin{IEEEproof}
The proof is deferred to the appendix.
\end{IEEEproof}
\end{lemma}

\section{Conclusion}
In conclusion, our study contributes to the ongoing discourse on joint probability mass function estimation by addressing certain aspects within the minimax framework, particularly for $l^p_p$ loss functions where $1 \le p \le 2$, including the important $l_1$ loss. Additionally, our findings extend the understanding of minimax optimality to encompass a range of $f$-divergences, such as KL, $\chi^2$, Squared Hellinger, and Le Cam divergences. By doing so, our work not only builds upon and extends existing methodologies in crucial directions.



\section*{Acknowledgment}
This work was supported in part by the ONR under Grant N00014-19-1-2621


%
\newpage
\bibliographystyle{IEEEtran}
\bibliography{references.bib}

\begin{thebibliography}{10}
\providecommand{\url}[1]{#1}
\csname url@samestyle\endcsname
\providecommand{\newblock}{\relax}
\providecommand{\bibinfo}[2]{#2}
\providecommand{\BIBentrySTDinterwordspacing}{\spaceskip=0pt\relax}
\providecommand{\BIBentryALTinterwordstretchfactor}{4}
\providecommand{\BIBentryALTinterwordspacing}{\spaceskip=\fontdimen2\font plus
\BIBentryALTinterwordstretchfactor\fontdimen3\font minus \fontdimen4\font\relax}
\providecommand{\BIBforeignlanguage}[2]{{%
\expandafter\ifx\csname l@#1\endcsname\relax
\typeout{** WARNING: IEEEtran.bst: No hyphenation pattern has been}%
\typeout{** loaded for the language `#1'. Using the pattern for}%
\typeout{** the default language instead.}%
\else
\language=\csname l@#1\endcsname
\fi
#2}}
\providecommand{\BIBdecl}{\relax}
\BIBdecl

\bibitem{Wald49}
\BIBentryALTinterwordspacing
A.~Wald, ``Statistical decision functions,'' \emph{The Annals of Mathematical Statistics}, vol.~20, no.~2, pp. 165--205, 1949. [Online]. Available: \url{http://www.jstor.org/stable/2236853}
\BIBentrySTDinterwordspacing

\bibitem{trybula1958some}
S.~Trybula, ``Some problems of simultaneous minimax estimation,'' \emph{The Annals of Mathematical Statistics}, vol.~29, no.~1, pp. 245--253, 1958.

\bibitem{Olkin1979}
\BIBentryALTinterwordspacing
I.~Olkin and M.~Sobel, ``Admissible and minimax estimation for the multinomial distribution and for k independent binomial distributions,'' \emph{The Annals of Statistics}, vol.~7, no.~2, Mar. 1979. [Online]. Available: \url{https://doi.org/10.1214/aos/1176344613}
\BIBentrySTDinterwordspacing

\bibitem{wilczynski1985minimax}
M.~Wilczy{\'n}ski, ``Minimax estimation for the multinomial and multivariate hypergeometric distributions,'' \emph{Sankhy{\=a}: The Indian Journal of Statistics, Series A}, pp. 128--132, 1985.

\bibitem{braess2004bernstein}
D.~Braess and T.~Sauer, ``Bernstein polynomials and learning theory,'' \emph{Journal of Approximation Theory}, vol. 128, no.~2, pp. 187--206, 2004.

\bibitem{HanJW14}
\BIBentryALTinterwordspacing
Y.~Han, J.~Jiao, and T.~Weissman, ``Minimax estimation of discrete distributions under {$l_1$} loss,'' \emph{CoRR}, vol. abs/1411.1467, 2014. [Online]. Available: \url{http://arxiv.org/abs/1411.1467}
\BIBentrySTDinterwordspacing

\bibitem{pmlr-v40-Kamath15}
\BIBentryALTinterwordspacing
S.~Kamath, A.~Orlitsky, D.~Pichapati, and A.~T. Suresh, ``On learning distributions from their samples,'' in \emph{Proceedings of The 28th Conference on Learning Theory}, ser. Proceedings of Machine Learning Research, P.~Grünwald, E.~Hazan, and S.~Kale, Eds., vol.~40.\hskip 1em plus 0.5em minus 0.4em\relax Paris, France: PMLR, 03--06 Jul 2015, pp. 1066--1100. [Online]. Available: \url{https://proceedings.mlr.press/v40/Kamath15.html}
\BIBentrySTDinterwordspacing

\bibitem{onsemsup}
H.~S. Melihcan~Erol, E.~Sula, and L.~Zheng, ``On semi-supervised estimation of distributions,'' in \emph{2023 {IEEE} International Symposium on Information Theory ({ISIT})}.\hskip 1em plus 0.5em minus 0.4em\relax IEEE, Jun. 2023.

\bibitem{devore1993constructive}
R.~A. DeVore and G.~G. Lorentz, \emph{Constructive approximation}.\hskip 1em plus 0.5em minus 0.4em\relax Springer Science \& Business Media, 1993, vol. 303.

\bibitem{Chung2002}
\BIBentryALTinterwordspacing
F.~Chung and L.~Lu, ``Connected components in random graphs with given expected degree sequences,'' \emph{Annals of Combinatorics}, vol.~6, no.~2, pp. 125--145, Nov. 2002. [Online]. Available: \url{https://doi.org/10.1007/pl00012580}
\BIBentrySTDinterwordspacing

\end{thebibliography}

\newpage\phantom{blabla}\newpage
\appendices
\section{Supplementary Results (Cont'd)}
\begin{thm}[\cite{devore1993constructive},Theorem-3.1]  \label{thm:approximation_bernstein} If $f$ is bounded on A, differentiable in some neighborhood of $x$, and has second derivative $f^{\prime \prime}(x)$ for some $x \in A$, then
$$
\lim _{n \rightarrow \infty} n\left[B_n(f, x)-f(x)\right]=\frac{x(1-x)}{2} f^{\prime \prime}(x) .
$$
\end{thm}

\begin{thm}[\cite{Chung2002}, lemma-2] \label{thm:binomial_tail}
 Let $X_1, \ldots, X_n$ be independent random variables with
$$
\operatorname{Pr}\left(X_i=1\right)=p_i, \quad \operatorname{Pr}\left(X_i=0\right)=1-p_i .
$$
We consider the sum $X=\sum_{i=1}^n X_i$, with expectation $\mathrm{E}(X)=\sum_{i=1}^n p_i$. Then, we have:
\begin{align*}
& \operatorname{Pr}(X \leq \mathrm{E}(X)-\lambda) \leq e^{-\lambda^2 / 2 \mathrm{E}(X)} \\
& \operatorname{Pr}(X \geq \mathrm{E}(X)+\lambda) \leq e^{-\frac{\lambda^2}{2(\mathbb{E}(X)+\lambda / 3)}}
\end{align*}
\end{thm} 
\section{Omitted Proofs}
\subsection{Proof for \lemref{lem:bernstein_trick}}
\begin{IEEEproof}
Fix $c > 0$, let $\delta_1,  \delta_2 > 0$,that will be determined later, we define $\Delta_p(x,i) \defeq \abs{x^p - \Big( \frac{i}{n} \Big)^{p}} $. As a result of triangular inequality:
\begin{IEEEeqnarray}{rCl}
  \abs{H^n_p(x) - G^n_p(x)} \le \sum^n_{i=0}  {n \choose i } x^i (1- x)^{n-i}  r^p_i\ \Delta_p(x,i) \label{eq:to_seperate} \IEEEeqnarraynumspace
\end{IEEEeqnarray}
Before analyzing this sum,  we note that by the mean value theorem, there exists $\xi \in (x \wedge \frac{i}{n}, x \vee \frac{i}{n} )$ hence $x^p - \left(\frac{i}{n} \right)^{p} = \left(x - \frac{i}{n}\right) \xi^{p-1}$, thus:
\begin{align} 
\Delta_p(x,i) = \Big|x^p - \Big(\frac{i}{n}\Big)^{p} \Big| \le p \Big|x -  \frac{i}{n} \Big|  \Big|x \vee \frac{i}{n} \Big|^{p-1} \label{eq:lem2:delta_upp}
\end{align}
Now we analyze the sum in \eqref{eq:to_seperate} over the ranges $ i \le nx$, $nx > i $ separately. 
 For the case $i \le nx$: 
    \begin{IEEEeqnarray}{rCl}
 &\ \ & \sum_{i\le nx }  {n \choose i } x^i (1- x)^{n-i}  r^p_i\ \Delta_p(x,i)  \IEEEnonumber   \\ 
    &=& \sum_{i \le n ( x - \delta_1)} {n \choose i } x^i (1-x)^{n-i} r^p_i \Delta_p(x, i) \IEEEnonumber   \\ 
 \  &\ +&  \sum_{n (x - \delta_1) < i < nx} {n \choose i } x^i (1-x)^{n-i} r^p_i \Delta_p(x, i) \label{eq:lem2:delta_subs} \\
 &\le& \sum_{i \le n  (x - \delta_1)} {n \choose i } x^i (1-x)^{n-i} 2  \IEEEnonumber \\ 
 &\ +& \sum_{n (x - \delta_1 ) \le i \le nx }  {n \choose i } x^i (1- x)^{n-i}  r^p_i \ p \abs{x - \frac{i}{n}} x^{p-1} \label{eq:lem2:delta_less_tail}
    \\ 
    &\le& 2 e^{-\frac{n}{x} \delta^2_1}  +
 \sum_{n (x - \delta_1 )\le i \le nx } p {n \choose i } x^i (1- x)^{n-i}  r^p_i \delta_1 x^{p-1}\quad \IEEEnonumber \\ 
  &=& 2 e^{-\frac{n}{x} \delta^2_1}    + \frac{ p\ \delta_1}{x} H^n_{p}(x) \label{eq:bound_game}
   \end{IEEEeqnarray}
In \eqref{eq:lem2:delta_subs} we use \eqref{eq:lem2:delta_upp}.  In \eqref{eq:lem2:delta_less_tail}, we bound the lower tail of the binomial via \thmref{thm:binomial_tail} and observe that $\abs{x - \frac{i}{n}} \le \delta_1$ in the range $n \cdot  (x  -\delta_1) \le i \le n \cdot x$. We also note that in \eqref{eq:lem2:delta_less_tail}, \  $r^p_i \le 2$ and $\Delta_p(x,i) \le 1$.
Now we choose $\delta_1 = c_1 \sqrt{- \frac{x}{n} \log \left( \frac{1}{x} H^n_p(x)  \right) }$ and obtain:
\begin{IEEEeqnarray}{rCl} 
\eqref{eq:bound_game} &\le& 2 \ \left(  \frac{H^n_p(x)}{x} \right)^{c^2_1} + p\ c_1 \sqrt{- \frac{1}{n x} \log \left( \frac{1}{x} H^n_p(x)  \right) } H^n_p(x) \IEEEeqnarraynumspace    \IEEEnonumber
\end{IEEEeqnarray}
Hence to establish the lemma, we first show that $\sqrt{-\frac{1}{nx} \log{\frac{1}{x} H^n_p(x) }} = O( \frac{1}{  \sqrt{\log{n}} } )$. 
\begin{IEEEeqnarray*}{rCl}
\frac{1}{x} H^n_p(x) &\simeq& C_p \sum^n_{i=0} \frac{x^{p-1}}{i^{\frac{p}{2}} + 1} {n \choose i} x^i (1-x)^{n- i} \\
&\ge& C_p \frac{1}{n^{\frac{p}{2}} + 1 }  \left(\frac{ c \log^2{n}}{n} \right)^{^{p-1}} \sum^n_{i=0} {n \choose i} x^i (1-x)^{n-i}  \\ \IEEEyesnumber \label{eq:lower_bound}\\
&=& C_p \frac{1}{n^{\frac{p}{2}} + 1 }  \left(\frac{ c \log^2{n}}{n} \right)^{^{p-1}}\ge \frac{k'}{n^{\frac{3}{2}p}} \IEEEyesnumber \label{eq:Hpnx_lower_bound}
\end{IEEEeqnarray*}
In \eqref{eq:lower_bound}, we notice $x \ge \frac{c \log^2{n}}{n}$. On the right-hand side of \eqref{eq:Hpnx_lower_bound} we collect the constants in $k'$.
Therefore we establish that:
\begin{IEEEeqnarray}{rCl}
\sqrt{-\frac{1}{n x } \log{ \frac{1}{x} H^n_p(x) } } \le \sqrt{k'' \frac{1}{n x } \log{n} } \le  \bigOh{ \frac{1}{\log(n)} } \quad \IEEEeqnarraynumspace \label{eq:resulting_delta}
\end{IEEEeqnarray}
where in the first inequality we use that $x \ge c \frac{\log{n}}{n}$ and collect the constants in $k''$. Secondly, we need to show that $\frac{H^n_p(x)}{x}$ decays sufficiently fast. To this end, we let $\delta'$ be a parameter to be determined later:
\begin{align}
\frac{H^n_p(x)}{x} &= \sum_{i}  \frac{C_p}{i^{\frac{p}{2}} + 1}  {n \choose i} x^{i + p - 1} (1-x)^{n -i } \nonumber \\ 
&\le C_p x^{p-1 } \Ed{I \sim \bin(n,x)}{\frac{1}{\sqrt{I} + 1} } \nonumber \\
&\le C_p x^{p-1} \left( \prob{ \abs{I - n x} \ge \delta' }  + \frac{1}{ \sqrt{nx - \delta' }  + 1} \right) \nonumber  \\ 
&\le C_p x^{p-1} \left( e^{-\frac{\delta'^2}{n}}  + \frac{1}{ \sqrt{nx - \delta' }  + 1} \right)  \nonumber \\
\intertext{hence setting $\delta' = x \cdot n^{0.75}$:}
&\le C_p x^{p-1} \left(  e^{-\sqrt{n} x^2 } + \frac{1}{\sqrt{n \cdot x}  + 1  } (1 + o(1))\right)  \nonumber \\
&\le C_p x^{p-1} \left(\frac{1}{\sqrt{n \cdot x}  + 1  } (1 + o(1))\right) \nonumber \\ 
&\le C_p \frac{1}{\sqrt{n}} x^{p-\frac{1}{2}} (1 + o(1)) \nonumber \\
&= \bigOh{ n^{-\frac{1}{2}} } \label{eq:idk:lastline}
\end{align}
\eqref{eq:idk:lastline} follows from that $p \ge 1$. Therefore by choosing $c_1$ large enough we ensure that $B \left(\frac{H^n_p(x)}{x}  \right)^{c^2_1} = o(n^{-\frac{p}{2}})$.
Whereas in the  second case $i > n x$:
\begin{IEEEeqnarray}{rCl}
\IEEEeqnarraymulticol{3}{l}{  \sum_{nx < i \le nx +  \delta_2 } {n \choose i} x^i (1 -x )^{n-i} r^p_i \Delta_p (x, i)  } \IEEEnonumber  \\  
&\quad+& \sum_{  i >  nx + \delta_2  } {n \choose i} x^{i} (1-x)^{n-i} r^p_i \Delta_p(x,i)  \\ 
&\le& \sum_{nx < i \le nx +  \delta_2 } {n \choose i} x^i (1 -x )^{n-i} \abs{x - \frac{i}{n}} \left( \frac{i}{n} \right)^{p-1}   \IEEEnonumber \\  
&\quad+&2  \sum_{  i >  nx + \delta_2  } {n \choose i} x^{i} (1-x)^{n-i}  \\ 
&\le& \delta_2 \sum_{nx < i \le nx +  \delta_2 } {n \choose i} x^i (1 -x )^{n-i} \left( x + \delta_2 \right)^{p-1}    + 2 e^{-\frac{n \delta^2_2}{2 (x + \frac{\delta_2}{3} )}} \IEEEnonumber \\ 
\IEEEeqnarraymulticol{3}{l}{\le \delta_2 2^{p-1} \sum_{nx < i \le nx +  \delta_2 } {n \choose i} x^i (1 -x )^{n-i}   \left( x \vee \delta_2 \right)^{p-1} + e^{-\frac{n \delta^2_2}{\frac{2}{3} (x \vee \delta_2 )}}  }  \IEEEnonumber \\
\IEEEeqnarraymulticol{3}{l}{  \le  \delta_2\  2^{p-1} \sum_{nx < i \le nx +  \delta_2 } {n \choose i} x^i (1 -x )^{n-i}   \left( x \right)^{p-1}  + 2 e^{-\frac{n \delta^2_2}{\frac{2}{3} x }}  } \IEEEyesnumber
\label{eq:thm2:lower:last}
 \end{IEEEeqnarray}
Each step is justified in the corresponding step in the analysis for the range $i \le nx$, except now we are using the upper tail in \thmref{thm:binomial_tail}. In \eqref{eq:thm2:lower:last}, we see that the problem is identical to \eqref{eq:bound_game} except for the constants. Hence we choose $\delta_2 = c_2 \sqrt{-\frac{x}{n} \log \left( \frac{1}{x} H^n_p(x)  \right) }$ and  by \eqref{eq:Hpnx_lower_bound} we have $\delta_2 \le  c_3 \sqrt{ x \frac{\log{n} }{n} }$. This ensures that $\delta_2 \le x $ when $x \ge c \frac{\log^2{n}}{n}$ and the step \eqref{eq:thm2:lower:last} is justified.
\end{IEEEproof}


\end{document}